\documentstyle[amssymb,amsfonts]{amsart}

\def\R{{\hbox{\bf R}}}

\font \roman = cmr10 at 10 true pt

\def\allt#1{%
\smash{
\vtop{%
     \ialign{%
        ##\crcr
        $\hfil\displaystyle{\tilde \forall}\hfil$\crcr%
        \noalign{\kern1.5pt\nointerlineskip}
        $\hfil\!\!#1\hfil$\crcr\noalign{\kern1.5pt}
        }
       }
      } \hbox{$\vphantom{#1}$}
     }

\def\be#1{\begin{equation}\label{#1}}
\def\bas{\begin{align*}}
\def\eas{\end{align*}}
\def\bi{\begin{itemize}}
\def\ei{\end{itemize}}

\def\dist{{\hbox{\roman dist}}}

\def\Z{{\hbox{\bf Z}}}
\def\eps{\varepsilon}
\newenvironment{proof}{\noindent {\bf Proof} }{\endprf\par}
\def \endprf{\hfill  {\vrule height6pt width6pt depth0pt}\medskip}
\def\emph#1{{\it #1}}
\def\textbf#1{{\bf #1}}

\theoremstyle{plain}
  \newtheorem{theorem}[subsection]{Theorem}
  \newtheorem{conjecture}[subsection]{Conjecture}

\theoremstyle{remark}

\theoremstyle{definition}
  \newtheorem{definition}[subsection]{Definition}

\include{psfig}

\begin{document}

\title[Recent progress on Kakeya]{Recent progress on the Kakeya conjecture}

\author{Nets Katz}
\address{Department of Mathematics, Washington University at St. Louis, St. Louis MO 63130}
\email{nets@@math.wustl.edu}

\author{Terence Tao}
\address{Department of Mathematics, UCLA, Los Angeles CA 90095-1555}
\email{tao@@math.ucla.edu}


\begin{abstract} 
We survey recent developments on the Kakeya problem.
\end{abstract}

\maketitle

\section{Introduction}

The purpose of this article is to survey the developments on the Kakeya problem in recent years, concentrating on the period after the excellent survey of Wolff \cite{wolff:survey}, and including some recent work by us in \cite{KT}, \cite{kt:kakeya2}.  The results covered here are discussed to some extent in the later surveys \cite{Bo}, \cite{wolff:icm}, \cite{tao:notices}.  We will focus on the standard Kakeya problem for line segments and not discuss other important variants (such as Kakeya estimates for circles, light rays, or $k$-planes; see e.g. \cite{wolff:survey}, \cite{wolff:cone}, \cite{wolff:smsub}, \cite{erdogan}).  We shall describe the various arguments in a rather informal fashion in order to not obscure the main ideas too much with technicalities.

For any $n \geq 2$, define an {\it Besicovitch set} to be a subset of $\R^n$ which contains a unit line segment in every direction.  An old construction of Besicovitch shows that such sets can have arbitrarily small measure in any dimension, and can even be made to be measure zero.  Intuitively, this states that it is possible to compress a large number of non-parallel unit line segments into an arbitrarily small set.  

The study of such sets originated in the Kakeya needle problem: what is the smallest amount of area needed to rotate a unit needle in the plane?  (The above construction shows that one can rotate a needle using arbitrarily small area).  Later work by Fefferman, C\'ordoba, Bourgain, and others have shown that these sets are connected to problems in oscillatory integrals, number theory, and PDE.  We shall not discuss these matters here, but refer the reader to \cite{tao:notices}.

In applications one wishes to obtain more quantitative understanding of this compression effect by introducing a spatial discretization.  For instance, one could replace unit line segments by $1 \times \delta$ tubes for some $0 < \delta \ll 1$ and ask for the optimal compression of these tubes.  Equivalently, one can ask for bounds of the volume of the $\delta$-neighbourhood of a Besicovitch set.

This problem can be phrased in terms of the {\it (upper) Minkowski dimension} of the Besicovitch set.  Recall that a bounded set $E$ has Minkowski dimension $\alpha$ or less if and only if for every $0 < \delta \ll 1$ and $0 < \eps \ll 1$, one can cover $E$ by at most $C_\eps \delta^{-\alpha+\eps}$ balls of radius $\delta$.  For any $1 \leq p \leq n$, let $K_M(p, n)$ denote the statement that all Besicovitch sets in $\R^n$ have Minkowski dimension at least $p$.

\begin{conjecture}[Kakeya conjecture, Minkowski version] We have $K_M(n, n)$.
\end{conjecture}

This conjecture is known to be true in two dimensions but only partial progress has been made in higher dimensions.

A slightly stronger version of this conjecture can be phrased by using the Hausdorff dimension. Recall that a bounded set $E$ has Hausdorff dimension $\alpha$ or less if and only if for every $0 < \delta \ll 1$ and $0 < \eps \ll 1$, one can cover $E$ by balls $\{B\}$ of radius $r(B) \leq \delta$ such that
$$ \sum_{B} r(B)^{\alpha-\eps} \leq C_\eps.$$
Clearly, the Hausdorff dimension is less than or equal to the Minkowski dimension.  Let $K_H(p, n)$ denote the statement that all Besicovitch sets in $\R^n$ have Hausdorff dimension at least $p$.

\begin{conjecture}[Kakeya conjecture, Hausdorff version] We have $K_H(n,n)$. 
\end{conjecture}

There is also a more quantitative version of these conjectures, which we shall write by discretizing to some scale $0 < \delta \ll 1$. Fix $\delta$, and let $\Omega$ be a maximal $\delta$-separated set of directions on the sphere $S^{n-1}$.  For each $\omega \in \Omega$, let $T_\omega$ be a $\delta \times 1$ tube in the direction $\omega$ in the ball.  Let $0 < \lambda \leq 1$, and let $\tilde T_\omega$ be a subset of $T_\omega$ such that $|\tilde T_\omega| = \lambda |T_\omega|$.  We use $X \lessapprox Y$ to denote the estimate $X \leq C_\eps \delta^{-\eps} Y$ for all $\eps > 0$.

We say that $K_X(p,n)$ holds if one has the estimate
\be{max-est}
| \bigcup_{\omega \in \Omega} \tilde T_\omega | \gtrapprox
\lambda^p \delta^{n-p}
\end{equation}
for all $0 < \delta \ll 1$, and choice of $T_\omega$, $\tilde T_\omega$.  The estimate $K_X(1,n)$ is trivial; the difficulty is in making $p$ large, as one needs to prevent the $\tilde T_\omega$ from overlapping too heavily.  Note that the estimate is sharp when $\lambda \sim \delta$ and the $\tilde T_\omega$ are $\delta$-balls centered at the origin.

\begin{conjecture}[Kakeya conjecture, maximal version] We have $K_X(p,n)$.
\end{conjecture}

These maximal function estimates can also be generalized to the slightly stronger x-ray estimates (see e.g. \cite{drury:xray}, \cite{cdrdef}, \cite{wolff:x-ray}, \cite{laba:x-ray}), but we shall not discuss these extensions here.

\begin{figure}\label{old}
\begin{tabular}{|l|l|l|l|} \hline
Result & Dimension & Value of $p$ & \\\hline
Minkowski/Hausdorff & $n=2$ & $2$ & Davies 1971 \cite{davies} \\
Maximal & $n=2$ & $2$ & Cordoba 1977 \cite{cordoba:sieve} \\
Minkowski/Hausdorff & $n \geq 3$ & $(n-1)/2 + 1$ & Drury 1983 \cite{drury:xray} \\
Maximal & $n \geq 3$ & $(n-1)/2 + 1$ & Christ-Duoandikoetxea-Rubio de Francia 1987 \cite{cdrdef}\\ \hline
All & $n \geq 3$ & $(n-1)/2 + 1 + \eps_n$ & Bourgain 1991 \cite{borg:kakeya}\\
All & $n \geq 3$ & $(n-2)/2 + 2$ & Wolff 1995 \cite{wolff:kakeya}\\ \hline
Minkowski/Hausdorff & $n > 26$ & $13(n-1)/25 + 1$ & Bourgain 1998 \cite{borg:high-dim} \\
Maximal & $n \gg 1$ & $(1+\eps)n/2$ & Bourgain 1998 \cite{borg:high-dim}\\
Hausdorff & $n > 12$ & $6(n-1)/11 + 1$ & Katz-Tao 1999 \cite{KT}\\
Minkowski & $n > 8$ & $4(n-1)/7 + 1$ & Katz-Tao 1999 \cite{KT}\\
\hline
\end{tabular}
\caption{Prior results on $K_X(p,n)$, $K_H(p,n)$, $K_M(p,n)$.}
\end{figure}

These estimates are related by the implications
$$ K_X(p,n) \implies K_H(p,n) \implies K_M(p,n).$$
Indeed, to obtain $K_H(p,n)$ one needs only to prove \eqref{max-est} for $\lambda \gtrsim 1 / \log(1/\delta)$, while to obtain $K_M(p,n)$ one needs only to prove \eqref{max-est} for $\lambda = 1$.  The former claim is shown by pigeonholing the balls covering the Besicovitch set depending on the dyadic scale of their radii.  (In fact, one only needs \eqref{max-est} for $\lambda \gtrsim 1/\log\log(1/\delta)$ provided that one is willing to concatenate a few dyadic scales together; see \cite{borg:high-dim}).  Thus the main distinction of the maximal version of the conjecture is that it also handles the case of small $\lambda$.

In two dimensions all three conjectures are known to be true; indeed, we now have very precise estimates on the exact dependence of the bounds on $\delta$. In three and higher dimensions the conjectures are still open; we summarize the current state of progress on these three estimates in Figure \ref{kfig}.  These results are far from sharp; in particular, there is a strong chance that the Hausdorff and maximal results can be improved to match the Minkowski ones.

\begin{figure}\label{kfig}
\begin{tabular}{|l|l|l|l|} \hline
Result & Dimension & Value of $p$ & \\\hline
Minkowski & $n=3$ & $5/2 + \eps$ & Katz-{\L}aba-Tao 1999 \cite{katzlabatao} \\
Minkowski & $n=4$ & $3 + \eps$ & {\L}aba-Tao 2000 \cite{laba:medium} \\
Minkowski & $23 \geq n \geq 5$ & $(2-\sqrt{2})(n-4)+3$ & Katz-Tao 2000 \cite{kt:kakeya2}\\
Minkowski & $n > 23$ & $(n-1)/\alpha + 1$ & Katz-Tao 2000 \cite{kt:kakeya2}\\ \hline
Hausdorff & $n = 3,4$ & $(n-2)/2 + 2$ & Wolff 1995 \cite{wolff:kakeya}\\
Hausdorff & $n \geq 5$ & $(2-\sqrt{2})(n-4)+3$ & Katz-Tao 2000 \cite{kt:kakeya2}\\ \hline
Maximal & $8 \geq n \geq 3$ & $(n-2)/2 + 2$ & Wolff 1995 \cite{wolff:kakeya}\\
Maximal & $n \geq 9$ & $4(n-1)/7 + 1$ & Katz-Tao 2000 \cite{kt:kakeya2}\\
\hline
\end{tabular}

\caption{Current state of progress on $K_X(p,n)$, $K_H(p,n)$, $K_M(p,n)$ in three and higher dimensions.  $\alpha = 1.675\ldots$ is the largest root of $\alpha^3 - 4\alpha + 2 = 0$.}
\end{figure}

As a historical note, the maximal version of the conjecture is so named because it implies a bound on the \emph{Kakeya maximal operator}
$$ f^*(\omega) = \sup_{l // \omega} \int_l f$$
where $\omega \in S^{n-1}$ and the supremum ranges over all lines parallel to $\omega$.  Indeed, $K_X(p,n)$ is equivalent to the estimate
$$ \| f^* \|_{L^p(S^{n-1})} \leq C_{\eps,p} \| f \|_{L^p_{n/p - 1 + \eps}(\R^n)}$$
for all $f$ supported on a bounded set; see e.g. \cite{borg:kakeya}.

\section{The geometric method}\label{geom-sec}

The known arguments for the Kakeya problem can be divided into three categories: geometric combinatorics arguments, arithmetic combinatorics arguments, and hybrids of the two.  We begin by reviewing the geometric arguments, which were developed earlier.

The first non-trivial case occurs in two dimensions.  In this case the conjectures are known to be true; the statements $K_M(2,2)$, $K_H(2,2)$ are due to Davies \cite{davies}, while $K_X(2,2)$ is due to C\'ordoba \cite{cordoba:sieve}.  The idea is to compute the quantity
$$ \| \sum_{\omega \in \Omega} \chi_{\tilde T_\omega} \|_2^2$$
directly, using the geometric observation
$$ |T_\omega \cap T_{\omega'}| \lesssim \frac{\delta^2}{\delta + \angle(\omega, \omega')}.$$
(This may be viewed as a quantitative version of the fact that two (transverse) lines only intersect in at most point).  The argument does not extend well to general dimensions, however it does give the useful estimate 
\be{2-plane}
| \bigcup_{\omega \in \Omega'} T_\omega| \approx \sum_{\omega \in \Omega'} |T_\omega| \sim \delta^{n-1} \# \Omega'
\end{equation}
whenever $T_\omega$ are a collection of $\delta$-tubes in $\R^n$ which lie in a 2-plane and which have a $\delta$-separated set of directions.  In other words, tubes in a 2-plane are essentially disjoint.

In higher dimensions, an argument of Drury \cite{drury:xray} gives $K_M((n+1)/2, n)$ and $K_H((n+1)/2, n)$; this was later extended by Christ, Duoandikoetxea, and Rubio de Francia \cite{cdrdef} to give $K_X((n+1)/2, n)$. the main geometric ingredient is the fact that any two (separated) points are connected by at most one line.  These estimates can be proved using the ``bush argument'' in \cite{borg:kakeya}, but we present an alternate ``slice argument'' from \cite{borg:high-dim} for $K_M((n+1)/2, n)$, which is the easiest of the three.  It suffices to show the estimate 
$$ | E | \gtrapprox \delta^{(n-1)/2}.$$
where $E$ is the set $\bigcup_{\omega \in \Omega} T_\omega$.

After some rescaling, we may assume that the tubes $T_\omega$ all intersect the hyperplanes $x_n = 0$ and $x_n = 1$ and make an angle of $\leq 1/10$ with the vertical $e_n$.  For all $0 \leq t \leq 1$, let $\tilde A_t$ denote the intersection of $E$ with $\{ x_n = t\}$; for generic values of $t$ we have $|\tilde A_t| \lesssim |E|$.

The sets $\tilde A_t$ are essentially unions of $n-1$-dimensional disks of radius $\delta$.  Thus, if we let $A_t$ be a maximal $\delta$-separated subset of $\tilde A_t$, then $\# A_t \sim \delta^{1-n} |\tilde A_t|$.  Thus, if we apply a generic rescaling, we may assume that $\# A_0, \# A_1 \lesssim \delta^{1-n} |E|$.  

On the other hand, each tube $T_\omega$ intersects $A_0$ in essentially one point, and similarly for $A_1$.  Thus each $T_\omega$ can be identified with a pair in $A_0 \times A_1$.  Because each pair of points in $A_0 \times A_1$ has essentially only one tube connecting it, we thus see that
$$ \delta^{1-n} \sim \# \Omega \leq \# (A_0 \times A_1) \lesssim (\delta^{1-n} |E|)^2$$
which is the desired estimate.

In 1991, Bourgain \cite{borg:kakeya} found a small improvement to the bound $K_X((n+1)/2, n)$ in all dimensions $n$.  This was then improved further in 1995 by Wolff \cite{wolff:kakeya}, who obtained $K_X((n+2)/2, n)$ for all $n$ (thus unifying these results with the two-dimensional theory).  Very roughly (and omitting all the technicalities), the idea is to modify the bush argument as follows.  We would like to show an estimate such as
$$ |E| \gtrapprox \delta^{(n-2)/2}.$$
Since $\sum_{\omega \in \Omega} \chi_{T_\omega}$ has an $L^1$ norm of $\sim 1$ and is supported on $E$, we see that every point in $E$ is contained in $\sim |E|^{-1}$ tubes on the average.

Now consider a ``hairbrush'', or more precisely the set of all tubes $T_\omega$ that pass through a given stem tube $T_{\omega_0}$.  By dividing $T_{\omega_0}$ up into about $\delta^{-1}$ balls of radius $\delta$ and applying the previous observation, we expect the hairbrush to consist of about $\delta^{-1} |E|^{-1}$ tubes.  Since all these tubes pass through $T_{\omega_0}$, they must essentially lie in some 2-plane containing $T_{\omega_0}$.  These 2-planes are mostly distinct, so we can subdivide the hairbrush into disjoint sub-collections of tubes, each of which lie in a single 2-plane.  Applying \eqref{2-plane} to each sub-collection and summing, we see that the volume of the hairbrush is at least
$\delta^{n-1} \delta^{-1} |E|^{-1}$.  On the other hand, the hairbrush must be contained in $E$.  Combining the two statements we obtain the result.

It appears difficult to improve the bound $(n+2)/2$ by geometric arguments, especially when $n=3$.  However, some refinements in other directions have been obtained by these types of arguments.  For instance, when $p=(n+2)/2$ then one can generalize \eqref{max-est} to
\be{xray-est}
| \bigcup_{\omega \in \Omega} \bigcup_{i=1}^m \tilde T_{\omega,i} | \gtrapprox m^{c_{n,p}} \lambda^p \delta^{n-p}
\end{equation}
where $c_{n,p} > 0$ is a constant depending on $n$, $p$, and for each $\omega$, the sets $\tilde T_{\omega,i}$ are density-$\lambda$ subsets of disjoint tubes $T_{\omega,i}$ oriented in the direction $\omega$.  See \cite{wolff:x-ray} for the $n=3$ case (with the sharp value $c_{3,5/2} = 1/4$) and \cite{laba:x-ray} for the general case.  This estimate is equivalent to a certain mixed-norm estimate for the x-ray transform.  Wolff \cite{wolff:x-ray} observed that the estimate \eqref{xray-est} implied a certain interesting refinement to $K_M(p)$.  Namely, if $E$ was a Besicovitch set of Minkowski dimension exactly $p$, and the $\delta$-neighbourhood of $E$ was essentially given $\bigcup_{\omega \in \Omega} T_\omega$, then the map $\omega \mapsto T_\omega$ was ``almost Lipschitz'' in the sense that for any $\delta < \rho < 1$ and any randomly chosen $\omega, \omega'$ with $|\omega - \omega'| \lesssim \rho$, one has $\dist(T_\omega, T_{\omega'}) \lessapprox \rho$ with probability $\approx 1$.  (The reason for this is that if this property failed, then one could pass to the $\rho$-neighbourhood of $E$ and use \eqref{xray-est} to show that the Minkowski dimension of $E$ had to be strictly larger than $p$).  This almost Lipschitz property is usually referred to as \emph{stickiness} (tubes which are nearly parallel must stick close to each other).

This stickiness property is especially useful at the scale $\rho := \sqrt{\delta}$.  (This intermediate scale seems to appear everywhere in this theory!).  Roughly speaking, stickiness then asserts that the tubes $T_\omega$ can be grouped into about $\rho^{1-n}$ groups of $\rho^{1-n}$ tubes each, such that each group of tubes is contained inside a $1 \times \rho$ tube.

By analyzing how these groups of tubes can intersect each other, one can derive further properties on these tubes.  For instance, in joint work with I. {\L}aba we have shown that for generic points $x$ in $E$, the set of $1 \times \rho$ tubes which contain $x$ are mostly contained in a $\rho$-neighbourhood of a hyperplane; this property has been dubbed ``planiness''.  Also, the intersection of $\bigcup_{\omega \in \Omega} \tilde T_\omega$ with any $\rho$-ball has a certain structure, namely that it is essentially the union of $\rho \times \rho \times \delta \times \ldots \times \delta$ slabs (this property we refer to as ``graininess'').  See \cite{katzlabatao} for a rigorous version of these statements in the $n=3$ case, and \cite{laba:medium} for the higher-dimensional case.  (Actually, in dimensions $n>3$ at least one of the planiness and graininess properties can be improved further (so that the planes have lower dimension, or the grains have higher dimension).  See \cite{laba:medium}). 

When $n > 4$ the properties of stickiness and planiness are strong enough to obtain a small improvement to Wolff's estimate in the Minkowski setting; in other words, we have $K_M((n+2)/2 + \eps_n, n)$ for some $\eps_n > 0$.  See \cite{laba:medium}.  In $n=4$ the same result obtains but one also needs to exploit the graininess property \cite{laba:medium}.  When $n=3$ the properties of stickiness, planiness, and graininess are not quite sufficient to obtain an improvement, and one must also introduce the arithmetic techniques of the next section.  See \cite{katzlabatao}.

We close this section by mentioning two fairly simple observations which have simplified some of the more technical issues in the field (although it is clear that these observations do not address the main geometric and arithmetic combinatorial issues).  The first observation is known as the ``two ends reduction'', and basically allows one to assume that the set $\tilde T_\omega$ is not concentrated in a small sub-tube of $T_\omega$, but is rather spread out throughout all of $T_\omega$.  The point is that if the sets $\tilde T_\omega$ were consistently concentrated in small sub-tubes, then one could pass to the scale of these sub-tubes and then rescale to obtain a better counter-example; see \cite{wolff:kakeya}.  The second observation is dual to the first, and is known as the ``bilinear reduction''.  Roughly, it allows one to assume that for most points $x \in \bigcup_{\omega \in \Omega} T_\omega$, the tubes $T_\omega$ that pass through $x$ are not concentrated in a small angular sector, but are spread out over all angles.  This is because if the tubes through generic points were consistently concentrated in narrow sectors, then one could pass to these sectors and then rescale to obtain a better counter-example; see \cite{tvv:bilinear}.

Another (equally informal) way of stating this is that the two ends reduction allows one to assume that generic points in a tube have separation $\approx 1$, while the bilinear reduction allows one to assume that generic tubes through a point have angular separation $\approx 1$.  One can push these observations a bit further in a non-rigorous fashion, and assert that the sets $\tilde T_\omega$ behave like self-similar fractals of dimension $\log(\lambda)/\log(\delta)$, while the sets $\{ \omega: T_\omega \ni x \}$ behave like self-similar fractals of dimension $\log |E| / \log \delta$.  This leads to some interesting questions concerning sets of Furstenburg type; see \cite{wolff:distance}, \cite{wolff:survey}.

\section{The arithmetic method}\label{arithmetic-sec}

The geometric methods such as those in Wolff \cite{wolff:kakeya} appear to be fairly efficient in low dimensions, but are not very satisfactory in very high dimensions.  In 1999 Bourgain \cite{borg:high-dim} introduced a new argument, based on the arithmetic combinatorics of sums and differences, which gave better results in high dimensions.  The connection between Kakeya problems and the combinatorics of addition can already be seen by considering the analogy between line segments and arithmetic progressions.  (Indeed, the Kakeya conjecture can be reformulated in terms of arithmetic progressions, and this can be used to connect the Kakeya conjecture to several difficult conjectures in number theory such as the Montgomery conjectures for generic Dirichlet series.  We will not discuss this connection here, but refer the reader to \cite{Bo}).

For simplicity we first begin by discussing the Minkowski dimension problem, and indicate the additional difficulties involved in the Hausdorff and maximal settings later in this section.

These arguments begin by analyzing the ``two-slice'' proof
of the estimate $K_M((n+1)/2, n)$ in more detail.  Recall that each tube $T_\omega$ was associated to a distinct element of the product set $A_0 \times A_1$.  Let $G$ denote the subset of $A_0 \times A_1$ generated by these tubes, thus 
\be{g-card}
\# G \sim \delta^{1-n}.
\end{equation}
Let $N$ denote the quantity $\delta^{1-n} |E|$, thus we have
\be{a-card}
\# A_t \lesssim N
\end{equation}
for generic values of $t$.  We shall gloss over the definition of ``generic'' and make the non-rigorous assumption that \eqref{a-card} in fact holds for all $t$.

We are interested in bounds of the form
\be{g-bound}
\# G \lessapprox N^\alpha,
\end{equation}
since this would then give
$$ \delta^{1-n} \lessapprox (\delta^{1-n} |E|)^\alpha$$
which eventually yields the estimate 
$$ K_M( \frac{n}{\alpha} + \frac{\alpha-1}{\alpha}, n ).$$
For instance, the bound $K_M((n+1)/2, n)$ comes from the trivial estimate 
\be{triv}
\# G \lesssim N^2.
\end{equation}
The aim is thus to make $\alpha$ as small as possible; an estimate with $\alpha = 1$ would solve the Kakeya conjecture, for Minkowski dimensions at least.

To improve upon the trivial bound \eqref{triv} we need to use some further properties of $G$.  Let $Z$ denote the vector space $\R^n$, thus $G \subset Z \times Z$.  Let $\pi_-: Z \times Z \to Z$ denote the subtraction map $\pi_-: (a,b) \mapsto a-b$.  From the fact that the tubes $T_\omega$ all point in different directions we see that
\be{one}
\pi_- \hbox{ is one-to-one on } G.
\end{equation}

For each $t \in [0,1]$, let $\pi_t: Z \times Z \to Z$ denote the map $\pi_t: (a,b) \to (1-t)a + tb$.  If $(a,b) \in G$, then there is a tube $T_\omega$ containing $a$ and $b$, and hence the intermediate point $(1-t)a + tb$.  Thus one essentially has
\be{map}
\pi_t \hbox{ maps } G \hbox{ to } A_t
\end{equation}
for all $t$.

In other words, the map $\pi_t$ compresses $G$ (which at present appears to have size about $N^2$) into the much smaller set $A_t$ (which has size about $N$).  For instance, applying this with $t = 1/2$ we obtain
$$ \# \{ a+b: (a,b) \in G \} \lesssim N.$$
Comparing this with \eqref{one}, we see that the differences of $G$ are all distinct, whereas the sums of $G$ have very large overlap.  Intuitively, one expects these facts to conflict with each other; this can already be seen from the elementary observation
observation
\be{sum-diff}
a + b = a' + b' \iff a - b' = a' - b.
\end{equation}

However, it is not entirely trivial to convert observations such as \eqref{sum-diff} into bounds of the form \eqref{g-bound}.  There is a substantial literature on the relative sizes of sum-sets $A_0 + A_1$ and difference sets $A_0 - A_1$ (see e.g. the excellent survey \cite{ruzsa}), but much less is known about \emph{partial} sum-sets and difference sets, when one only considers a subset $G$ of pairs $A_0 \times A_1$.  One result in this direction is the Balog-Szemer\'edi theorem \cite{balog}, which asserts that if a large subset of $A_0 \times A_1$ has a small sum-set, then there exists large subsets $A'_0$, $A'_1$ of $A_0$, $A_1$ respectively such that $A'_0 + A'_1$ is also small.  In principle, this theorem should be quite useful for us; however, the quantitative bounds given by \cite{balog} were far too poor (the constants blow up much faster than exponential) to give any significant improvement to the estimate on $|E|$.
	
The main breakthrough came from Gowers \cite{gowers}, who developed a quantitative version of the Balog-Szemer\'edi theorem while working on the apparently unrelated problem of locating arithmetic progressions of length 4.  This argument was then adapted by Bourgain \cite{borg:high-dim} for the Kakeya problem.  

To state Bourgain's main ``sums-differences'' estimate, we pause to give some notation.

\begin{definition}  Let $t_1, \ldots, t_r$ be a sequence of reals and $1 \leq \alpha \leq 2$.  We say that the estimate $SD(t_1, \ldots, t_r; \alpha)$ holds if one has the bound \eqref{g-bound} for all integers $N > 0$, all vector spaces $Z$, all finite sets $G \subset Z \times Z$, and all finite subsets $A_{t_1}, \ldots, A_{t_r}$ of $Z$ such that \eqref{a-card}, \eqref{one}, \eqref{map} hold for $t = t_1, \ldots, t_r$.  We say that the estimate $SD(\alpha)$ holds if for every $\eps > 0$ there exists $t_{1,\eps}, \ldots, t_{r_\eps,\eps}$ such that $SD(t_{1,\eps}, \ldots, t_{r_\eps,\eps}; \alpha + \eps)$ holds.
\end{definition}

We thus have
$$ SD(\alpha) \implies K_M( \frac{n}{\alpha} + \frac{\alpha-1}{\alpha}, n )$$
for all $n$.  The trivial estimate \eqref{triv} becomes $SD(0,1; 2)$ in this notation.

By adapting Gowers' arguments, Bourgain showed $SD(0,1/2,1; 2 - \frac{1}{13})$, which thus implies $K_M( \frac{13 n + 12}{25}, n)$; note that this improves upon Wolff's bound $K_M( (n+2)/2, n)$ for $n > 26$.  The idea is to show that for many elements $(a,b) \in G$, there are many ways to represent $a-b$ as
$$ a-b = (a_1 - b_1) - (a_2 - b_2) + (a_3 - b_3)$$
where $a_1, a_2, a_3 \in A_0$ and $b_1, b_2, b_3 \in A_1$.  This limits the total number of possible values of $a-b$, which then limits the size of $G$.  The representations are obtained via the identity
$$ a-b = (a - b') - (a' - b') + (a' - b)$$
and the fact (from \eqref{sum-diff}) that there are generically a lot of solutions to the equations $a-b' = a_1 - b_1$, $a'-b' = a_2 - b_2$, $a' - b = a_3 - b_3$.  See \cite{borg:high-dim}; an alternate version (with $1/13$ replaced by an unspecified epsilon) is in \cite{katzlabatao}.

In \cite{KT} we developed a somewhat different approach to proving results of the form $SD(t_1, \ldots, t_r; \alpha)$.  The idea is to identify certain configurations of points in $G$ (e.g. trapezia, vertical line segments, corners, etc.) and obtain both lower and upper bounds for the number of such configurations in terms of $\# G$ and $N$.  Comparing these bounds then gives an estimate of the form \eqref{g-bound}.  Of course, one needs to utilize \eqref{one} in order to obtain an improvement over \eqref{triv}; this shall be accomplished by means of elementary linear algebra, re-expressing $\pi_-$ in terms of the other projections $\pi_t$.

We illustrate this technique with 

\begin{theorem}\label{kt-thm}\cite{KT} We have $SD(0,1/2,2/3,1; 2 - \frac{1}{4})$.  In particular, we have $K_M((4n+3)/7, n)$.  
\end{theorem}

\begin{proof}
Let $Z$, $N$, $G$, $A_0$, $A_{1/2}$, $A_{2/3}$, $A_1$ be as above.
Define a \emph{vertical line segment} to be a pair $(g_1, g_2) \in G \times G$
such that $\pi_0(g_1) = \pi_0(g_2)$.  (The notation comes from depicting $A_0$, $A_1$ as one-dimensional sets, so that $G \subset A_0 \times A_1$ becomes a two-dimensional set and $+_0$ is the projection to the first co-ordinate).  Let $V$ denote the space of all vertical line segments. From an easy Cauchy-Schwarz argument we have
$$ \# V \geq (\# G)^2 / \# A_0 \gtrsim (\# G)^2 / N.$$

Now define a \emph{trapezoid} to be a pair $( (g_1, g_2), (g_3, g_4) ) \in V \times V$ of vertical line segments such that $\pi_1(g_1) = \pi_1(g_3)$ and $\pi_{2/3}(g_2) = \pi_{2/3}(g_4)$, and let $T$ denote the space of all trapezoids.  From Cauchy-Schwarz we have a lower bound for the cardinality of $T$:
\be{t-lower}
\# T \geq (\# V)^2 / (\# A_1 \# A_{2/3}) \gtrsim (\# V)^2 / N^2 \gtrsim (\# G)^4 / N^4.
\end{equation}
On the other hand, we can obtain an upper bound for the cardinality of $T$ as follows.  Consider the map $f: T \to A_{1/2} \times A_{1/2} \times A_1$ defined by
$$ f( (g_1, g_2), (g_3, g_4) ) := ( \pi_{1/2}(g_1), \pi_{1/2}(g_2), \pi_1(g_4) ).$$
We claim that this map is one-to-one.  At first glance, this seems unlikely, as $T$ has four degrees of freedom with respect to $Z$ ($G$ has 2 degrees of freedom, hence $V$ has $2*2-1 = 3$, hence $T$ has $2*3-2 = 4$) and $f$ only specifies three of these four degrees.  However, we can use \eqref{one} to recover the fourth degree of freedom.  Specifically, we take advantage of the identity
\be{ident}
\pi_-(g_3) = -2 \pi_{1/2}(g_1) + 4 \pi_{1/2}(g_2) - 2 \pi_1(g_4)
\end{equation}
for all $((g_1,g_2),(g_3,g_4)) \in T$, to conclude that the quantity $f((g_1,g_2),(g_3,g_4))$ determines $\pi_-(g_3)$, which then determines $g_3$ by \eqref{one}.  From some straightforward linear algebra one can then check that $f((g_1,g_2), (g_3, g_4))$ determines all of $((g_1,g_2), (g_3, g_4))$, or in other words that $f$ is one-to-one.  Hence
$$ \# T \leq \# A_{1/2} \# A_{1/2} \# A_{2/3} \lesssim N^3.$$
Combining this with \eqref{t-lower} we obtain \eqref{g-bound} with $\alpha = 2-1/4$ as claimed.

\end{proof}

A similar argument gives $SD(0, 1/2, 1; 2 - \frac{1}{6})$; see \cite{KT}.

We now sketch the more recent results in \cite{kt:kakeya2}.  The starting point is the observation that one can improve upon the estimate $SD(1.75)$ in Theorem \ref{kt-thm} by analyzing the proof more carefully.  Define the map $\nu: V \to Z$ by
\be{nu-1}
\nu( g_1, g_2 ) := -2 \pi_{1/2}(g_1) + 4 \pi_{1/2}(g_2).
\end{equation}
Since $\pi_0(g_1) = \pi_0(g_2)$, $\nu$ can also be written as
\be{nu-2}
\nu( g_1, g_2 ) = - \pi_1(g_1) + 3 \pi_{2/3}(g_2)
\end{equation}
or
\be{nu-3}
\nu( g_1, g_2 ) = \pi_-(g_1) + 2 \pi_1(g_2)
\end{equation}

From \eqref{nu-2} we see that 
$$ \nu(g_1, g_2) = \nu(g_3, g_4)$$
whenever $(g_1, g_2), (g_3, g_4)) \in T$.  From this and \eqref{nu-1}, \eqref{nu-2}, \eqref{nu-3} we obtain a derivation of \eqref{ident}.

In light of this, it is natural to consider for each $\nu_0 \in Z$ the set $$V_{\nu_0} := \{ (g_1, g_2): \nu(g_1, g_2) = \nu_0 \}.$$
From any of \eqref{nu-1}, \eqref{nu-2}, \eqref{nu-3} we see that $(g_1, g_2)$ is completely determined by $g_2$ and $\nu_0$.  Thus if we set
$$ \tilde V_{\nu_0} := \{ g_2: (g_1, g_2) \in V_{\nu_0} \}$$
then there is a one-to-one correspondence between $\tilde V_{\nu_0}$ and $V_{\nu_0}$.

On the other hand, we see from the previous discussion on $T$ that
$$ \# \{ ((g_1, g_2), (g_3, g_4)) \in V^2: \nu(g_1,g_2) = \nu(g_3,g_4);
\pi_{2/3}(g_2) = \pi_{2/3}(g_4) \} $$
$$\geq \# T \gtrsim \# V \frac{\# V}{N^2} \gtrsim \# V (\# G)^2 N^{-3}.$$
We can re-arrange this as
$$ \sum_{\nu_0 \in Z} \# \{ ((g_1, g_2), (g_3, g_4)) \in V_{\nu_0}^2:
\pi_{2/3}(g_2) = \pi_{2/3}(g_4) \} \gtrsim \sum_{\nu_0 \in \Z} \# V_{\nu_0} (\# G)^2 N^{-3},$$
which can be re-arranged further as
$$ \sum_{\nu_0 \in Z} \# \{ (g_2, g_4) \in \tilde V_{\nu_0}^2:
\pi_{2/3}(g_2) = \pi_{2/3}(g_4) \} \gtrsim \sum_{\nu_0 \in \Z} \# \tilde V_{\nu_0} (\# G)^2 N^{-3}.$$
Thus for generic choices of $\nu_0$, we expect
$$
\# \{ (g_2, g_4) \in \tilde V_{\nu_0}^2: \pi_{2/3}(g_2) = \pi_{2/3}(g_4) \} \gtrsim \# \tilde V_{\nu_0} (\# G)^2 N^{-3},$$
which morally speaking should imply
$$ \# \{ \pi_{2/3}(g_2): g_2 \in \tilde V_{\nu_0} \} \lesssim \frac{\# \tilde V_{\nu_0}}{ (\# G)^2 N^{-3} }.$$
(Strictly speaking, we need to pass to a large subset of $\# \tilde V_{\nu_0}$ to achieve this, but we gloss over this technicality).  Thus if we fix $\nu_0$ and define
$$ G^* := \tilde V_{\nu_0}$$
$$ A^*_t := \{ \pi_t(g): g \in G \}$$
and
\be{n-def}
N^* := \frac{\# G^*}{ (\# G)^2 N^{-3} }
\end{equation}
then we have just shown that
$$ \# A^*_{2/3} \lesssim N^*.$$
A similar argument using \eqref{nu-1} gives
$$ \# A^*_{1/2} \lesssim N^*$$
(at least if $\nu_0$ is chosen sufficiently generically).  In fact, we can show
\be{as-card}
\# A^*_t \lesssim N^*
\end{equation}
for all $t \neq 1$ by this method, providing that \eqref{a-card} holds for the correct values of $t$.

On the other hand, from \eqref{nu-3} and \eqref{one} we see that for $(g_1, g_2) \in V$, $\pi_1(g_2)$ and $\nu(g_1, g_2)$ determine $g_1$, and hence all of $(g_1, g_2)$.  Thus we have
\be{vnu-bound}
\# G^* = \# V_{\nu_0} \leq \# A_1 \lesssim N.
\end{equation}

Observe by construction that \eqref{one} and \eqref{map} hold with $G$, $A_t$ replaced by $G^*$, $A^*_t$.  (Indeed, the former is clear since $G^* \subset G$).  Thus if we have a suitable bound $SD(\alpha)$, we may apply it to $G^*$, $N^*$, and conclude
$$ \# G^* \lesssim (\# N^*)^{\alpha}.$$
Combining this with \eqref{n-def} and \eqref{vnu-bound} we eventually obtain
$$ \# G \lesssim N^{2 - 1/(2\alpha)}.$$
One can make this rigorous and conclude that

\begin{theorem}\label{imply}\cite{kt:kakeya2} $SD(\alpha)$ implies $SD(2 - \frac{1}{2\alpha})$.
\end{theorem}

Thus $SD(2)$ gives $SD(2 - \frac{1}{4})$, which gives $SD(2 - \frac{2}{7})$, and so forth.  Iterating this and taking limits we obtain $SD(1 + \sqrt{2}/2) = SD(1.707\ldots)$, which implies
$$ K_M((2-\sqrt{2})n + (\sqrt{2}-1)).$$

We have been able to set up a more sophisticated iteration scheme along the same lines, where the role of the vertical line segments $(g_1, g_2)$ are replaced by ``corners'' $(g_1, g_2, g_3)$, where $\pi_0(g_1) = \pi_0(g_2)$ and $\pi_1(g_2) = \pi_1(g_3)$.  By similar arguments to the above, we have been able to obtain $SD(\alpha)$, where $\alpha = 1.675\ldots$ solves the cubic $\alpha^3 - 4\alpha + 2 = 0$.  In high dimensions ($n \geq 24$) this gives the best progress on the Minkowski version of the Kakeya version to date.  (For the $n < 24$ results, see below).

We have certainly not exhausted all the possibilities of this approach, and it is quite likely that one can improve the above results by finding some good upper and lower bounds for various configurations of objects in $G$.  To date, the arithmetic techniques used so far have all been remarkably elementary (with some minor combinatorial technicalities in making such terms as ``generic'' rigorous).

We now discuss the issue of how well these Minkowski arguments transfer over to the more difficult Hausdorff or maximal Kakeya problems.  In these cases the tubes $T_\omega$ are replaced by subsets $\tilde T_\omega$ of density $\lambda$.  
Let $E$ denote the union of the $\tilde T_\omega$, and let $A_t$ be defined by $E$ as before.  We can still ensure an estimate of the form \eqref{a-card}, however a difficulty arises in obtaining \eqref{map} for $t \neq 0,1$.  This is because statements $a, b \in \tilde T_\omega$ no longer imply that $(1-t)a + tb \in \tilde T_\omega$.

In the special case $t=1/2$ this amounts to the problem of locating arithmetic progressions of length 3 in $\tilde T_\omega$.  A famous theorem of Roth asserts that this is possible providing the density $\lambda$ is sufficiently large; an argument of Heath-Brown \cite{heath-brown} actually shows that one can take $\lambda \gtrsim 1/(\log 1/\delta)^c$ for some small absolute constant $c$.  This turns out to be enough to handle the Hausdorff problem (see \cite{borg:high-dim}), but is far too weak to give maximal function results\footnote{In \cite{borg:high-dim} Bourgain was still able to achieve a small maximal function improvement, namely $K_X((1/2 + \eps)n, n)$ for some absolute constant $\eps > 0$, by locating triples whose \emph{reciprocals} were in arithmetic progression (which are more common, due to the curvature of $1/x$), and then running similar arguments to the above.}.  Unfortunately these ideas do not seem to easily give Hausdorff or maximal results when one needs several values of $t$ at once, because of the difficulty of finding arithmetic progressions of length 4 or higher (see e.g. \cite{gowers}).

Nevertheless, one can adapt the proof of $SD(2 - 1/4)$ or $SD(1 + \sqrt{2}/2)$ to give Hausdorff and maximal results.  (It is likely that one can also do this for the more complicated argument of $SD(1.675\ldots)$, but we have not checked this).  The key observation is that one has a substantial amount of freedom to choose which slices $A_t$ to use in these arguments.  For instance, the argument giving $SD(0,1,1/2,2/3;2 - 1/4)$ can easily be modified to give
$SD(0,1,t,\frac{1}{2-t}; 2-1/4)$ for any $0 < t < 1$.  One can then show that for generic values of $t$, the statements \eqref{map} are mostly true for $t$ and $1/(2-t)$ (but with some loss depending on $\lambda$), which is enough to run the argument properly.  (For small $\lambda$ it is convenient to take advantage of the ``two-ends'' reduction mentioned in the previous section, to ensure that one can still find good slices $A_0$, $A_1$ of roughly unit separation).

Note that none of these results can reproduce Wolff's result in low dimensions.  This seems to be due to the fact that these arithmetic arguments do not fully exploit the Besicovitch set, but only use a finite number of slices of it.  In high dimensions this is less of an issue since each slice only has one less dimension than the full set, but this loss of information becomes very significant in two and three dimensions.  In the next section we mention some recent attempts to recapture this loss.

\section{Hybrid methods}\label{hybrid-sec}

As we have seen, there are two known ways of improving over the benchmark exponent of $(n+1)/2$.  Geometric techniques such as those of Wolff can improve the constant term to $(n+2)/2$, while arithmetic techniques can improve the linear term, obtaining an exponent of the form $n/\alpha + (\alpha-1)/\alpha$.  It is natural to ask whether one can combine the two techniques and obtain new Kakeya estimates which improve upon both the linear and the constant term.

One such attempt, joint with I. {\L}aba, is in \cite{katzlabatao}, in which the properties of stickiness, planiness, and graininess are used to obtain additional structure on the slices $A_t$.  However, this argument has so far only proved useful for the three-dimensional Minkowski problem.

Another approach, which we develop in \cite{kt:kakeya2}, is to ``de-slice'' the arguments of the previous section, and try to exploit the Besicovitch set more fully.  This turns out to be fairly straightforward.  For instance, one can give a re-proof of

\begin{theorem}\label{deslice}\cite{KT} We have $K_M((4n+3) / 7, n)$
\end{theorem}

\begin{proof} (Informal)
Let $E$ be the union of all the tubes $T_\omega$ as before; this set can be thought of as the union of about $\delta^{-n}|E|$ $\delta$-balls.  Since each tube occupies about $\delta^{-1}$ of these balls, we can expect two random $T$, $T'$ to have a probability of about $\delta^{n-2}/|E|$ of intersecting.

Define a \emph{quadrilateral} to be four tubes $T_1$, $T_2$, $T_3$, $T_4$ in $E$ such that $T_i$ intersects $T_{i+1}$ in a $\delta$-ball for $i=1,2,3,4$ (with $T_5 := T_1$). The total number of 4-tuples of tubes is $\delta^{4(1-n)}$, but there are four constraints, and so by the above heuristic we expect the total number of quadrilaterals is bounded below by $\delta^{4(1-n)} (\delta^{n-2}/|E|)^4$.  (This can be proven by two applications of Cauchy-Schwarz as before).  Let $x_i$ be the center of the ball where $T_i$ and $T_{i+1}$ intersect.

A quadrilateral is mostly determined by the three points (which are usually in $E$):
$$  \frac{1}{2} x_1 + \frac{1}{2} x_2, 2x_2 - x_3, \frac{2}{3} x_3 + \frac{1}{3} x_4.$$
Given these three points, one can determine $x_1 - x_4$ as a linear combination.  This determines the tube $T_4$, because the tubes point in different directions.  Thus $x_4$ (for instance) has only one degree of freedom, and once this degree is specified one can reconstruct the entire quadrilateral.  Thus the number of quadrilaterals is bounded by $(\delta^{-n}|E|)^3 \delta^{-1}$.

Combining these bounds gives the same bound $(4n+3)/7$ as previously.
\end{proof}

This argument should be compared with the proof of Theorem \ref{kt-thm}.

There are several advantages to this de-sliced formulation.  Firstly, it becomes somewhat easier to tackle the Hausdorff and maximal versions of the Kakeya conjecture in this setting, as one has less of a need to locate arithmetic progressions.  Secondly, one is now able to utilize such facts as the C\'ordoba estimate \eqref{2-plane} (which is very difficult to use in the purely arithmetic setting).  Indeed, recall from our discussion of Wolff's theorem that the tubes which pass through a single stem tube are essentially disjoint.  Another way of saying this is that given any stem tube $T$ and any point $x$ at a distance $\sim 1$ from $T$, there should only be $O(1)$ tubes which intersect both $T$ and $x$.  This fact can be used to improve slightly the bounds in the proof of Theorem \ref{deslice}, and can eventually lead to the improvement of $K_H((4n+5)/7, n)$.  A similar argument (but combined with a lifting argument from $n$ dimensions to $n+1$ dimensions, in order to fully exploit C\'ordoba's observation) can be used to convert the sliced result $SD(1 + \sqrt{2}/2)$ to the estimate 
$$ K_H( (2 - \sqrt{2}) (n-4) + 3, n),$$
which is currently the best estimate known on the Kakeya problem in dimensions $5 \leq n \leq 23$.

\end{document}